\documentclass[12pt]{amsart}
\usepackage{amsthm}
\usepackage{tikz-cd}
\usepackage{amssymb}
\usepackage{enumerate}
\usepackage{amsmath} 
\usepackage[mathscr]{euscript}
\makeatletter
\@namedef{subjclassname@2010}{%
  \textup{2010} Mathematics Subject Classification}
\makeatother
\newtheorem{thm}{Theorem}[section]
\newtheorem{lem}[thm]{Lemma}
\newtheorem{prop}[thm]{Proposition}
\newtheorem{defi}[thm]{Definition}

\newtheorem{corl}[thm]{Corollary}
\newtheorem{xrem}[thm]{Remark}
\newtheorem{exm}[thm]{Example}
\begin{document}

\baselineskip=17pt

\author{Nabanita Ray}
\address{Institute of Mathematical Sciences\\ CIT Campus, Taramani, Chennai 600113, India and Homi Bhabha National Institute, Training School Complex, Anushakti Nagar, Mumbai 400094, 
India}
\email[Nabanita Ray]{nabanitar@imsc.res.in}

\begin{abstract}
In this paper, we prove that  blown up at seven general points admits a conic
bundle structure over $\mathbb{P}^1$and it can be embedded as $(2, 2)$ divisor in $\mathbb{P}^1\times\mathbb{P}^2$. Conversely, any
smooth surface in the complete linear system $\mid (2, 2) \mid$ of $\mathbb{P}^1\times\mathbb{P}^2$ can be obtained as an embedding of blowing up
$\mathbb{P}^ 2$ at seven points. We also show that smooth surface linearly equivalent to $(2, 2)$ in $\mathbb{P}^1\times\mathbb{P}^2$
has at most four $(-2)$ curves .

 \end{abstract}

\title{Geometry of $\mathbb{P}^{2}$ blown up at seven points}

\maketitle
\vskip 4mm

\section{Introduction}

It is well known that $\mathbb{P}^{2}$ blown up at six general points is isomorphic to a smooth cubic in
$\mathbb{P}^{3}$ and the embedding is given by the anti-canonical divisor. Conversely, any smooth cubic
of $\mathbb{P}^{3}$ is isomorphic to $\mathbb{P}^{2}$ blown up at six general points. Similarly, it is interesting to ask if
$\mathbb{P}^{2}$ blown up at seven general points can be embedded in a three-fold. It has been shown in \cite{BA} that $4\pi^*H - 2E_i -\sum^{ 7}_{ j=1,j\neq i} E_j$ is a very ample divisor of $\mathbb{P}^{2}$ blown up at seven general
points. We will show in this note that this very ample divisor gives an embedding of $\mathbb{P}^{2}$
blown up at seven general points in $\mathbb{P}^{1}\times\mathbb{P}^{2}$ as well as in $\mathbb{P}^{5}$ . We will also see that $\pi^* H - E_i$ gives the first projection to $\mathbb{P}^{1}$ and the anti-canonical divisor gives the second projection to
$\mathbb{P}^{2}$ .

Griffiths and Harris (\cite{GH} page 546) have shown that the map given by the anti-canonical divisor from $\mathbb{P}^{2}$ blown up at seven general points to $\mathbb{P}^{2}$ gives a double cover of $\mathbb{P}^{2}$ . We show that the
anti-canonical is the only linear system which expresses $\mathbb{P}^{2}$ blown up at seven general points
as a 2-sheeted branch cover of $\mathbb{P}^{2}$ .

We consider here the conic bundle structure over a variety, as defined e.g. in \cite{S}: a triple $(V, S, \pi)$,
where $\pi : V \rightarrow S$ is a rational map, whose generic fiber is an irreducible rational curve and
$S$ is a nonsingular variety, is called a conic bundle over the base $S$ or simply a conic over $S$.
In our case the morphism from $\mathbb{P}^{2}$ blown up at seven general points admits a conic bundle
structure over  $\mathbb{P}^{1}$ given by the complete linear system $\pi^* H - E_i$ . Here, we list all the linear
systems which gives a conic bundle structure of  $\mathbb{P}^{2}$ blown up at seven general points over
 $\mathbb{P}^{1}$ (See Theorem \ref{thm4.1}). We also see that if we embedded  $\mathbb{P}^{2}$ blown up at seven general points in
$\mathbb{P}^{1}\times\mathbb{P}^{2}$ , then lines of the surface are components of non irreducible fibers corresponding to
the first projection i.e. the conic bundle map.(See Section [6]).

In \cite{S} Sarkisov shows, if $(V, S, \pi)$ is a conic bundle and $\mathscr{E} = \pi_* \mathscr{O}_V (-K_V )$ is a locally free
sheaf of rank three on $S$, then $V$ can be embedded into $\mathbb{P}(\mathscr{E})$ and each fiber of the morphism $\pi$ is rational curve of degree two in $\mathbb{P}(\mathscr{E})$. Any rational curve of degree two is a
conic in some $\mathbb{P}^{2}$ (\cite{H}.IV.3.). In our case $\mathscr{E} = \mathscr{O}_{\mathbb{P}^1}\oplus\mathscr{O}_{\mathbb{P}^1} \oplus\mathscr{O}_{\mathbb{P}^1}$ then $\mathbb{P}(\mathscr{E}) = \mathbb{P}^{1}\times\mathbb{P}^{2}$ and the embedding of
$\mathbb{P}^{2}$ blown up at seven “general points” is linearly equivalent to a $(2, 2)$ divisor in $\mathbb{P}^{1}\times\mathbb{P}^{2}$ since
Pic($\mathbb{P}^{1}\times\mathbb{P}^{2}$ ) =$\mathbb{Z}\oplus\mathbb{Z}$. Conversely, we will see that any smooth surface of $\mid (2, 2) \mid$ is isomorphic
to  $\mathbb{P}^{2}$  blown up at seven points. Also there are examples of smooth surfaces linearly equivalent
to $(2, 2)$ of $\mathbb{P}^{1}\times\mathbb{P}^{2}$ which are isomorphic to $\mathbb{P}^{2}$ blown up at seven non-general points.

Here, we will see that any smooth surface $S \sim (2, 2)$ of $\mathbb{P}^{1}\times\mathbb{P}^{2}$ and let $ C$ be a curve in
$S$, then $C.C \geq-2$. Also, we show that there will be at most four curves in $S$ which have
self-intersection $(-2)$.

\subsection*{ Acknowledgement}
I would like to thank my advisor D.S Nagaraj for his valuable guidance and constructive suggestions throughout this project. This work is financially supported  by a fellowship from IMSc,Chennai (HBNI), DAE, Government of India.

\section{Notations and Definitions}
We denote by $\mathbb{P}^2$ the projective plane over the field $\mathbb{C}$ of complex numbers
and by $\widetilde{\mathbb{P}^2_r} $ or $\widetilde{\mathbb{P}}^2_{P_1P_2...P_r} $ the blow up of $\mathbb{P}^2$ at $r$ points $P_1, P_2,...,P_r$. Let
$\pi$:$\widetilde{\mathbb{P}^2_r}$ $\rightarrow$ $\mathbb{ P}^2$ be the blow up map, H be a hyperplane section of $\mathbb{P}^2$ and $E_i$ be the exceptional curve 
corresponding to the point $P_i$, for $i=1,2,...,r$. The \textit{ Picard group} of $\widetilde{\mathbb{P}^2_r} $ is $\mathbb{Z}.l\oplus_{i=1} ^r \mathbb{Z}. E_i$, where $l$ is $\pi^\ast (H)$.

Let $P_1,P_2,...,P_r \in \mathbb{P}^2$. An \textit{ admissible transformation} is a quadratic transformation of $\mathbb{P}^2$ centered at  three of the $P_i$'s (call them $P_1,P_2,P_3$) 
 . This gives a new $\mathbb{P}^2$ and a new sets of $r$ points, 
namely $Q_1,Q_2,Q_3$ and the image of $P_4,P_5,...,P_r$. We say $P_1, P_2,...,P_r$ are in \textit{ general position} if no three are collinear and furthermore, after any finite sequence of 
admissible transformations, the new set of $r$ points also has no three collinear. See (\cite{H}.V.5)  for more details.

\begin{defi}
 A triple $(V, S, \pi)$, where $\pi: V \rightarrow S$ is a rational map whose generic fiber is an irreducible rational curve and $S$ is a nonsingular variety, is called a \textit{
conic bundle over the base $S$} or simply a  \textit{ conic over $S$}.
\end{defi}

\begin{defi}
 A conic $(V, S, \pi)$ is called  regular  if the map $\pi: V \rightarrow S$ is a flat morphism of nonsingular varieties. A regular conic is called  standard  if 
the morphism $\pi: V \rightarrow S$ is relatively minimal, i.e. $\pi^{-1}(D)\subset V$ is an irreducible divisor for each irreducible divisor $D\subset S$.

\end{defi}
For more details on conic bundles see \cite{S}. In this note we define honest conic bundle.
\begin{defi}
 A regular conic is called  honest  if $\pi^{-1}(D)\subset V$ is a reduced  divisor for each integral divisor $D\subset S$.
\end{defi}

\begin{xrem}
 If $(V, S, \pi)$ is a honest conic over $S$ and $S\simeq\mathbb{P}^1$, then $D\sim \{pt\}$ and each fiber is a reduced conic.
\end{xrem}

\section{$\mathbb{P}^2$ blown up at seven general points as a double cover of $\mathbb{P}^2$}

In Section 3, 4, 5 and 6 we denote by $\widetilde{\mathbb{P}^2_r} $ or $\widetilde{\mathbb{P}}^2_{P_1P_2...P_r} $, $\mathbb{P}^2$ blown up at $r$ \it general points \rm otherwise it will be stated  where
$P_1,P_2,...,P_r\in \mathbb{P}^2$.

\begin{lem}\label{lem3.1}
 Let $\widetilde{\mathbb{ P}^2_7} $ be $\mathbb{P}^2$ blow-up at seven general points. If there is a degree two finite map from 
 $\widetilde{\mathbb{ P}^2_7} $ to $\mathbb{P}^2$, then this map can only be given by the anti-canonical 
 divisor, which is $\mid3 \pi^\ast H -\sum_{i=1}^{7}E_i\mid$, where $H$ is a hyperplane section of $\mathbb{P}^2$.

\begin{proof}
 First, we claim that any finite degree two map from $\widetilde{\mathbb{ P}^2_7} $ to $\mathbb{P}^2$ is defined by a complete linear system. Assume that the claim
is not true. Let $\mathfrak{d}$ be a sub-linear system of the complete linear system $\mid D \mid $, which defines a finite degree two map 
$i_{\mathfrak{d}}:\widetilde{\mathbb{ P}^2_7} \rightarrow\mathbb{P}^2$. As $\mathfrak{d}$ is base point free, $\mid D \mid $ is also base point
free. Hence $\mid D \mid $ induces the map $i_{\mid D \mid}:\widetilde{\mathbb{ P}^2_7} \rightarrow\mathbb{P}^N$, where $h^0(\widetilde{\mathbb{ P}^2_7},\mid D \mid)=N+1$. 
Note that, $2\geq$ dim$(i_{\mid D \mid}(\widetilde{\mathbb{ P}^2_7}))\geq$ dim$(i_{\mathfrak{d}}(\widetilde{\mathbb{ P}^2_7}))=2$. This implies $i_{\mid D \mid}(\widetilde{\mathbb{ P}
^2_7})$ is a non-degenerated surface in $\mathbb{P}^N$, $i_{\mid D \mid}$ is generically finite map and $i_{\mid D \mid}^\ast\mathcal{O}_{\mathbb{P}^N}(1)=\mathcal{O}_{\widetilde{\mathbb{ P}^2_7}}( D)$. As 
$\mathfrak{d}\subseteq \mid D\mid$ corresponds to the degree two map $i_{\mathfrak{d}}$,  $D^2=2$. Hence, 
\begin{center}
$D^2=
i_{\mid D \mid}^\ast\mathcal{O}_{\mathbb{P}^N}(1)\cdot i_{\mid D \mid}^\ast\mathcal{O}_{\mathbb{P}^N}(1)$,\\
$ D^2=deg (i_{\mid D \mid})\bigr(\mathcal{O}_{\mathbb{P}^N}(1)\mid_{i_{\mid D \mid}
(\widetilde{\mathbb{ P}^2_7})}\cdot\mathcal{O}_{\mathbb{P}^N}(1)\mid_{i_{\mid D \mid}(\widetilde{\mathbb{ P}^2_7})}\bigl)$,
\\$2=deg (i_{\mid D \mid})\cdot deg(i_{\mid D \mid}(\widetilde{\mathbb{ P}^2_7}))$.
\end{center}

Hence, either $deg (i_{\mid D \mid})=2$, or $deg (i_{\mid D \mid})=1$.

If $deg (i_{\mid D \mid})=2$, then $deg (i_{\mid D \mid}(\widetilde{\mathbb{ P}^2_7}))=1$ in $\mathbb{P}^N$. 
If $Y\subset\mathbb{P}^N $ is an irreducible non-degenerated surface of degree $d$, then $2+d-1\geq N$( See \rm[GH] page 173). In our 
case $d=1$, hence $N=2$. Therefore, dim$\mathfrak{d} = $ dim$\mid D \mid$, and hence the linear system $\mathfrak{d}$ is a complete linear system.

If $deg (i_{\mid D \mid})=1$, then $deg (i_{\mid D \mid}(\widetilde{\mathbb{ P}^2_7}))=2$ in $\mathbb{P}^N$. Using the same  result as referred in the
above paragraph, 
we get $N\leq 3$. The case $N=2$ is not possible, because there is no degree 2 surface in $\mathbb{P}^2$. Now if $N=3$, then 
$i_{\mid D \mid}:\widetilde{\mathbb{ P}^2_7} \rightarrow\mathbb{P}^3$ is a degree one map and the image is degree two surface in
$\mathbb{P}^3$. Up to isomorphism,
there are two irreducible degree two surfaces in $\mathbb{P}^3$, one are smooth quadrics which are isomorphic to $\mathbb{P}^1\times\mathbb{P}^1$ and 
the other are
cone over plane conic curve which have a singularity at the vertex. Let $P\in \mathbb{P}^3\backslash i_{\mid D \mid}(\widetilde{\mathbb{ P}^2_7})$ and
we take a projection from $P$ to a
hyperplane in $\mathbb{P}^3$, $p:\mathbb{P}^3\backslash\{P\}\rightarrow\mathbb{P}^2$ such that we have $i_{\mathfrak{d}}=p\circ i_{\mid D \mid}$, where
$\widetilde{\mathbb{ P}^2_7}\xrightarrow{i_{\mid D \mid}}i_{\mid D \mid}(\widetilde{\mathbb{ P}^2_7})\xrightarrow{p}\mathbb{P}^2$. As $i_{\mathfrak{d}}$ is
a finite map of degree two, 
$i_{\mid D \mid}$ is a finite map of degree one and $p$ is a finite degree two map. But there is not any finite degree one map from
$\widetilde{\mathbb{ P}^2_7}$ to either $\mathbb{P}^1\times\mathbb{P}^1$ or from  cone over plane conic curve. Hence we conclude that, 
there doesn't exist any such degree one map $\widetilde{\mathbb{ P}^2_7}\xrightarrow{i_{\mid D \mid}} \mathbb{P}^3$. Our claim is proved i.e., 
any finite degree two map from $\widetilde{\mathbb{ P}^2_7} $ to $\mathbb{P}^2$ 
is determined by a complete linear system.

Now, let $\phi : \widetilde{\mathbb{ P}^2_7} \rightarrow \mathbb{ P}^2$ be a degree two map defined by the complete linear system $\mid D \mid $ 
where $D= a \pi^\ast H -\sum_{i=1}^{7} b_i E_i$
, $a>0$, $b_i\geq 0 $ and at least one $b_i > 0$. As $\phi$ is a degree two map, $D^2=2$ which implies $a^2-\sum_{i=1}^{7} b_i^2 = 2$.  
We know that expected dimension, expdim$\mid D\mid \leq$ dim$\mid D \mid=2$.\\ 
expdim$\mid D\mid= \frac{(a+1)(a+2)}{2}-1-\sum_{i=1}^{7}\frac{b_i(b_i+1)}{2}\leq2$\\
$\Rightarrow a^2+3a+2-2-\sum b_i^2- \sum b_i\leq 4$\\
$\Rightarrow 3a+2-\sum b_i\leq 4$ (as $a^2-\sum_{i=1}^{7} b_i^2 = 2 $)\\ 
$\Rightarrow 3a-2\leq \sum b_i$

If $\{x_i\}$ and $\{y_i\}$ are two real sequences, then by the Schwarz's inequality,\\ $\mid\sum_{i}x_iy_i\leq\mid\sum_{i}x_i^2\mid .\mid\sum_{i}y_i^2\mid$.\\
We replace $x_i=1, y_i=b_i$ for $i=1,2,\cdots,7$ and $x_i=0, y_i=0$ for $i>$7. Then we have\\
$(\sum b_i)^2 \leq 7(\sum b_i^2)$\\
$\Rightarrow (3a-2)^2\leq7(a^2-2)$\\
$\Rightarrow a^2-6a+9\leq0$\\
$\Rightarrow (a-3)^2\leq0 $\\
$\Rightarrow a=3 $.

Hence the only possibility of $D$ is $3 \pi^\ast H -\sum_{i=1}^{7}E_i$. So the   degree two map $\phi$ is given by the divisor 
$D=3 \pi^\ast H -\sum_{i=1}^{7}E_i$.  This gives the proof of the statement of the theorem.

\end{proof}

\end{lem}

\section{conic bundle structure of $\widetilde{\mathbb{ P}^2_7} $ over $\mathbb{P}^1$}

\begin{thm}\label{thm4.1}
If $\widetilde{\mathbb{ P}^2_7} $ admits a conic bundle structure over $\mathbb{ P}^1$ given by the linear system $\mid D\mid$. 
Then D will have one of the following forms,

(I) $\pi^\ast H-E_i$, $1\leq i \leq 7$,

(II) 2$\pi^\ast H-\sum_{i=1}^{4}E_{l_i}$, $l_i$ are distinct and $1\leq l_i \leq 7$, 

(III) 3$\pi^\ast H-2 E_i-\sum_{j=2}^{6}E_{k_j}$, $i$ and $k_j$ are distinct $1\leq i,k_j \leq 7$, 

(IV) $4\pi^\ast H-\sum_{j=1}^{4}E_{k_j}-\sum_{i=5}^{7}2E_{l_i}$ where $k_j$ and $l_i$ are distinct and $1\leq l_i, k_j \leq 7$ and 

(V) 5$\pi^\ast H- E_i-\sum_{j=2}^{7}2E_{k_j}$ where $i$ and $k_j$ are distinct and $1\leq i,k_j \leq 7$.
 \begin{proof}
We claim that any conic bundle map from $\widetilde{\mathbb{ P}^2_7}$ to $\mathbb{ P}^1$ is defined by a complete linear system. Assume that the claim
is not true. Let $\mathfrak{b}$ be a sub-linear system of the complete linear system $\mid D\mid$ which gives a conic bundle map 
$j_{\mathfrak{b}}:\widetilde{\mathbb{ P}^2_7}\rightarrow\mathbb{P}^1$. Now, we define a morphism $j_{\mid D \mid}:\widetilde{\mathbb{ P}^2_7} \rightarrow\mathbb{P}^n$ where 
$h^0(\widetilde{\mathbb{ P}^2_7},\mid D \mid)=n+1$. Here we have $2 \geq $ dim$(j_{\mid D \mid}(\widetilde{\mathbb{ P}^2_7})) \geq $ dim$(j_{\mathfrak{b}}(\widetilde{\mathbb{ P}^2_7}))=1$.

If dim$j_{\mid D \mid}(\widetilde{\mathbb{ P}^2_7})=2$, then $j_{\mid D \mid}$ is a generically finite map between two surfaces,
hence $D^2= deg(j_{\mid D \mid}).deg(j_{\mid D \mid}(\widetilde{\mathbb{ P}^2_7}))$. Also we have $deg(j_{\mid D \mid})>0$ and $deg(j_{\mid D \mid}(\widetilde{\mathbb{ P}^2_7}))>0$.
But we have $j_{\mathfrak{b}}^\ast\mathscr{O}_{\mathbb{P}^1}(1)=D$ and $j_{\mathfrak{b}}$ is a conic bundle map hence the fibers of $j_{\mathfrak{b}}$  are disjoint curves linearly 
equivalent to $D$. Hence $D^2=0$ which gives a contradiction.
 
So dim$j_{\mid D \mid}(\widetilde{\mathbb{ P}^2_7})=1$. We can project repeatedly from outside the image of $j_{\mid D \mid}$ and can get the following commutative diagram.

\[
\begin{tikzcd}
\widetilde{\mathbb{ P}^2_7} \arrow{r}{j_{\mid D \mid}} \arrow{dr}{j_{\mathfrak{b}}}
& j_{\mid D \mid}(\widetilde{\mathbb{ P}^2_7}) \arrow{d}{q}\\
& \mathbb{P}^{1}
\end{tikzcd}
\]

where $q$ is a finite map between two curves. But ${j_{\mathfrak{b}}}$ has generically connected fiber so, deg$(q)$=1.

Now, let $\psi : \widetilde{\mathbb{ P}^2_7} \rightarrow \psi ( \widetilde{\mathbb{ P}^2_7})\subseteq\mathbb{ P}^n$ be a conic bundle map over a rational curve corresponding to the 
divisor $D$. So a general element of $\mid D\mid$ is a smooth conic. Genus of the curve =$ g(D)=0$. $\psi^\ast(x)\sim D$. Any two fibre do not intersect each other. So, $D^2$=0 
and dim$\mid D\mid$=$n\geq1$. Let $ D= a \pi^\ast H -\sum_{i=1}^{7} b_i E_i$ where $a>$0
and $b_i\geq0$. 

$D^2$=0  $\Rightarrow a^2= \sum b_i^2 $

$g(D)=0$

$\Rightarrow\frac{(a-1)(a-2)}{2}-\sum_{i=1}^{7}\frac{b_i(b_i-1)}{2}$=0 

$\Rightarrow a^2-3a+2-\sum b_i^2 + \sum b_i = 0$

$\Rightarrow 3a-\sum b_i=2$

By Schwarz's inequality $\mid\sum_{i}x_iy_i\leq\mid\sum_{i}x_i^2\mid .\mid\sum_{i}y_i^2\mid$

Here we replace $x_i=1, y_i=b_i$ for $i=1,2,...,7 $ and $x_i=0, y_i=0$ for $i>7$. So we have

$\Rightarrow(\sum b_i)^2 \leq 7(\sum b_i^2)$

$\Rightarrow (3a-2)^2\leq7a^2$

$\Rightarrow a^2-6a+2\leq0$

$\Rightarrow a<6$

So possible values of $a$ are 1,2,3,4,5.

 Case I ($a=1$)

 If $a$=1 then $b_i$=1 and $b_j$=0 where $i\neq j$. So w.l.o.g take $b_1=$1, then $D=\pi^\ast H-E_1$. First, we need to check $\mid D\mid$ gives a map to $\mathbb{P}^1$.
1=expdim$\mid D \mid\leq$ dim($ \mid D \mid$). Now we claim dim($ \mid D \mid$)=1. But this is clear because, any curve of $\mid D\mid$ corresponds to a line passing through $P_1$ in
$\mathbb{P}^2$. So $\mid D\mid$ gives a map to $\mathbb{P}^1$ and generic fiber is an irreducible rational curve.
Hence ($\widetilde{\mathbb{ P}^2_7} , \mathbb{P}^1, \pi^\ast H-E_i$) is conic over $\mathbb{P}^1$.

Case II ($a=2$)

 If $a=2$, only possibilities of $b_i$ are $b_{k_{j}}$=1 where $j=1,2,3,4$ and others are zero. Then w.l.o.g we can consider 
$D=2\pi^\ast H-\sum_{i=1}^{4}E_i$. But dim$(\mid D \mid)=1$ because any curve of $\mid D\mid$ corresponds to a conic in $\mathbb{P}^2$ passing through $P_1, P_2, P_3, P_4$ which are in
general position. Hence, $D=2\pi^\ast H-\sum_{i=1}^{4}E_i$ gives a map to $\mathbb{P}^1$ where generic fiber is an irreducible rational curve. Then ($\widetilde{\mathbb{ P}^2_7} , 
\mathbb{P}^1, 2\pi^\ast H-\sum_{i=1}^{4}E_i$) is also conic bundle over $\mathbb{P}^1$.

Case III($a=3$)

If $a=3$, then only possibilities of $b_i$'s are $b_{k_{j}}$=1 where $j=1,2,3,4,5$ and among the other two  $b_i$ one is two and the other is zero. 
So w.l.o.g we consider $D=3\pi^\ast H-2 E_1-\sum_{i=2}^{6}E_i$ and by similar argument dim$\mid D\mid$=1.

So this $D$ gives a map to $\mathbb{P}^1$ where a generic fiber is an irreducible rational curve. Then ($\widetilde{\mathbb{ P}^2_7} , 
\mathbb{P}^1, 3\pi^\ast H-2 E_1-\sum_{i=2}^{6}E_i$) is also conic bundle over $\mathbb{P}^1$.

Case IV($a=4$)

 If $a=4$, then we have $10=\sum_{i=1}^7 b_i$ and $16=\sum_{i=1}^7 b_i^2$. Only possibilities of $b_i$'s are $b_{i_{1}}=b_{i_{2}}=b_{i_{3}}=b_{i_{4}}$=1 and 
$b_{j_{5}}=b_{j_{6}}=b_{j_{7}}$=2, where $i_k$ and $ j_l$ are distinct. So w.l.o.g consider $D=4\pi^\ast H-\sum_{j=1}^{4}E_j-\sum_{i=5}^{7}2E_i$. Now take a quadratic transformation $\phi$
of $\mathbb{P}^2$ centered at $P_5,P_6,P_7$ is $\phi : \mathbb{P}^2 \dashrightarrow \mathbb{P}^2$. Then $\phi (P_i)=P_i'$ for $i=1,2,3,4$. 
Here $\pi$ is the blow up map at the points $P_1,P_2,...,P_7$ of $\mathbb{P}^2$ and $L_{ij}$ is strict transformation of the line joining $P_i$ and $P_j$ in 
$\widetilde{\mathbb{P}^2}_{P_1P_2...P_7}$. $\pi'$ is the blow up map at
the points $P_1',P_2',P_3',P_4',Q_5,Q_6,Q_7$ of $\mathbb{P}^2$ where $\pi'(L_{56})=Q_7,\pi'(L_{67})=Q_5, 
\pi'(L_{57})=Q_6$, then $\widetilde{\mathbb{P}^2}_{P_1P_2...P_7}=\widetilde{\mathbb{P}^2}_{P_1'P_2'P_3'P_4'Q_5Q_6Q_7}$. $D \sim D'$ where 
$D'=2\pi'^\ast H'-\sum_{i=1}^{4}E_i'$  (\cite{H}.V.4.). Hence $\mid D\mid=\mid D'\mid$ and we have proved in case II, $\mid D'\mid$ gives map to $\mathbb{P}^1$ and  generic fiber 
is irreducible rational curve.

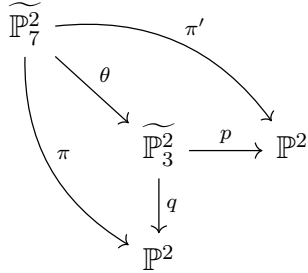
\[\begin{tikzcd}
\widetilde{\mathbb{P}^2_7}
\arrow[drr, bend left, "\pi'"]
\arrow[ddr, bend right, "\pi"]
\arrow[dr,"\theta"] & & \\
& \widetilde{\mathbb{P}^2_3} \arrow[r, "p"] \arrow[d, "q"]
& \mathbb{P}^2  \\
& \mathbb{P}^2 
\end{tikzcd}\]

Then similarly ($\widetilde{\mathbb{ P}^2_7} , \mathbb{P}^1,4\pi^\ast H-\sum_{j=1}^{4}E_j-\sum_{i=5}^{7}2E_i $) is also conic bundle over $\mathbb{P}^1$

Case V($a=5$)

If $a=5$, then we have $13=\sum_{i=1}^7 b_i$ and $25=\sum_{i=1}^7 b_i^2$. Only possibilities of $b_i$'s are $b_i=1, b_{j_{k}}=2$ where $k=1,2,...,6$ and $i\neq j_{k}$. In particular 
$D=5\pi^\ast H- E_1-\sum_{i=2}^{7}2E_i$. Now we take a quadratic transformation of $\mathbb{P}^2$ centered at $P_2,P_3,P_4$. After quadratic transformation 
we get $\mid D'\mid=\mid D \mid$ where $D'=4\pi'^\ast H'-\sum_{j=1}^{4}E_j'-\sum_{i=5}^{7}2E_i'$. Then using the  Case IV we have that ($\widetilde{\mathbb{ P}^2_7} , \mathbb{P}^1,
5\pi^\ast H- E_1-\sum_{i=2}^{7}2E_i$) is a conic bundle over $\mathbb{P}^1$.

Hence the result follows.
 \end{proof}

\end{thm}

\section{$\widetilde{\mathbb{P}^2_7} $ embedded as a (2,2) divisor of $\mathbb{P}^{1}\times\mathbb{P}^{2}$}
\begin{thm}\label{thm5.1}
 $\mathbb{P}^{2}$ blown up at seven general points can be embedded as a (2,2) divisor of $\mathbb{P}^{1}\times\mathbb{P}^{2}$.
 \begin{proof}
  We have the morphism $p_1:\widetilde{\mathbb{P}^2_7}\rightarrow\mathbb{P}^1$ which is defined by the linear system $\mid\pi^\ast H-E_1\mid$ and the morphism 
$p_2:\widetilde{\mathbb{P}^2_7}$ $\rightarrow$ $\mathbb{P}^2$ which is defined by the linear system $\mid 3 \pi^\ast H -\sum_{i=1}^{7}E_i\mid$.

Then we can define a morphism $p_1\times p_2: \widetilde{\mathbb{P}^2_7} \longrightarrow \mathbb{P}^{1}\times\mathbb{P}^{2}$. We have the Segre embedding, $\nu: \mathbb{P}^{1}\times\mathbb{P}^{2}\hookrightarrow \mathbb{P}^5 $ by the very ample divisor $\mid (1,1)\mid$. Hence the morphism $\nu\circ (p_1\times p_2): \widetilde{\mathbb{P}^2_7}\rightarrow\mathbb{P}^5$ is given by the linear system $\mid 4 \pi^\ast H - 2 E_1-\sum_{i=2}^{7}E_i\mid$

But we know $4 \pi^\ast H - 2 E_1-\sum_{i=2}^{7}E_i$ is a very ample divisor of $\widetilde{\mathbb{P}^2_7} $ by (Theorem 2.1 (\cite{Ha}).
Therefore $\nu\circ (p_1\times p_2)$ gives an embedding and $p_1\times p_2: \widetilde{\mathbb{P}^2_7} \longrightarrow \mathbb{P}^{1}\times\mathbb{P}^{2}$ is a closed immersion. 

With a slight abuse of notations, we will still use $ \widetilde{\mathbb{P}^2_7}$ and $\mathbb{P}^{1}\times\mathbb{P}^{2}$ to denote their embeddings into $\mathbb{P}^5$.

So $\widetilde{\mathbb{P}^2_7}$ is a non-singular
surface in a threefold and it corresponds to an element of its Weil divisor group. 
We know Pic($\mathbb{P}^{1}\times\mathbb{P}^{2}$)= Pic$(\mathbb{P}^{1})\oplus$Pic$(\mathbb{P}^{2})= \mathbb{Z}\oplus\mathbb{ Z}$
and the generators of the group are pt$\times\mathbb{P}^{2}$=(1,0) and $\mathbb{P}^1\times$H=(0,1).
Let $\widetilde{\mathbb{P}^2_7}\sim(a,b)$ in Pic($\mathbb{P}^{1}\times\mathbb{P}^{2}$).  $(a,b)=(a.pt\times\mathbb{P}^2) +( \mathbb{P}^1\times  bH)$. So generic fiber of the
first projection from $(a,b)$ is a curve of degree $b$ and generic fiber of the second projection contains $a$ number of points. 
So from Lemma(\ref{lem3.1}) and Theorem(\ref{thm4.1}) we have $a=2$ and $b=2$.
Hence the result is proved.
 \end{proof}

\end{thm}

\begin{xrem}
 Similarly the divisors $ 4 \pi^\ast H - 2 E_i-\sum_{j=1,j\neq i}^{7}E_j$ also give an embedding of $\widetilde{\mathbb{P}^2_7}$ in
 $\mathbb{P}^{1}\times\mathbb{P}^{2}$ as a smooth surface of $(2,2)$ type.
\end{xrem}
\begin{xrem}\label{xrem5.3}
 Let $(\widetilde{\mathbb{ P}^2}_{P_1P_2...P_7}, \mathbb{P}^1, f)$  be a conic over $\mathbb{P}^1$ where the morphism $f$ is defined by the divisor $\pi ^\ast H-E_1$. Each fiber
 is $f^*(P)\sim\pi ^\ast H-E_1$, where $P\sim\mathscr{O}_{\mathbb{P}^1}(1)$. 
 As $(\pi ^\ast H-E_1).(4 \pi^\ast H - 2 E_1-\sum_{i=2}^{7}E_i)=2$ , each fiber of $f$ is a degree two rational curve in $\mathbb{P}^5$ i.e a conic in some plane of  $\mathbb{P}^5$.
\end{xrem}

In the Theorem (\ref{thm5.1}) we have seen that the divisor of Lemma (\ref{lem3.1}) along with divisors of Case-I of Theorem (\ref{thm4.1}) give us an embedding.
Now we are interested to know  whether divisors of the other 
cases of Theorem (\ref{thm4.1}) along with the unique degree two map of Lemma (\ref{lem3.1}) will give us an embedding in $\mathbb{P}^{1}\times\mathbb{P}^{2}$ or not.

\begin{xrem}\label{xrem5.4}
Case I (7 possibilities)

\[\begin{tikzcd}
\widetilde{\mathbb{ P}^2}_{P_1P_2...P_7}
\arrow[r, "\mid D'\mid"] 
\arrow[rd, "\mid D''\mid"]
& \mathbb{P}^1 \\
& \mathbb{P}^2
\end{tikzcd}\]

$D'=\pi ^\ast H-E_i$ and $D''=3 \pi^\ast H -\sum_{i=1}^{7}E_i$

Theorem (\ref{thm5.1}) implies $D'+D''$ gives an closed immersion.

\subparagraph*{}
Case II (35 possibilities)

$D'=2\pi^\ast H-\sum_{i=1}^{4}E_{P_{j_i}}$ and $D''=3 \pi^\ast H -\sum_{i=1}^{7}E_{P_i}$. In particular consider $D'=2\pi^\ast H-\sum_{i=1}^{4}E_{P_i}$. 
Now take a quadratic transformation centered at $P_1,P_2,P_3$ and we will get $\widetilde{\mathbb{P}^2}_{P_1P_2...P_7}=\widetilde{\mathbb{P}^2}_{Q_1Q_2Q_3P_4'P_5'P_6'P_7'}$ 
and $D'\sim F'$ and $D''\sim F''$ where $F'=\pi'^\ast H'-E'_{P_4'}$
and $F''=3 \pi'^\ast H' -\sum_{i=1}^{3}E'_{Q_i}-\sum_{i=4}^{7}E'_{P_i'}$ and $\pi':\widetilde{\mathbb{P}^2}_{Q_1Q_2Q_3P_4'P_5'P_6'P_7'}\longrightarrow \mathbb{P}^2$ is the blowing up map
at the points $Q_1,Q_2,Q_3,P_4',P_5',P_6',P_7'$. From Case I we know, $F'+F''$ is very ample divisor hence $D'+D''$ also a very ample divisor. 
So $\widetilde{\mathbb{ P}^2}_{P_1P_2...P_7}\xrightarrow{D'+D''}\mathbb{P}^{1}\times\mathbb{P}^{2}$ is an closed immersion.

\subparagraph*{} 
Case III (42 possibilities)

$D'=3\pi^\ast H-2 E_i-\sum_{j=1,k_j\neq i}^{5}E_{k_j}$ and $D''=3 \pi^\ast H -\sum_{i=1}^{7}E_{P_i}$, in particular $D'=3\pi^\ast H-2 E_{P_1}-\sum_{j=2}^{6}E_{P_j}$. Now take the quadratic 
transformation centered at $P_1,P_2,P_3$,  we will get $\widetilde{\mathbb{P}^2}_{P_1P_2...P_7}=\widetilde{\mathbb{P}^2}_{Q_1Q_2Q_3P_4'P_5'P_6'P_7'}$ and $D'\sim F'$ and 
$D''\sim F''$, where $F'=2\pi'^\ast H'-E'_{Q_1}-E'_{P_4'}-E'_{P_5'}-E'_{P_6'}$ and $F''=3 \pi'^\ast H' -\sum_{i=1}^{3}E'_{Q_i}-\sum_{i=4}^{7}E'_{P_i'}$. 
Now, we are in the same position as Case II and repeating the Case II, there exist $ R_1,R_2,...R_7\in \mathbb{P}^2$ such that $\pi''^\ast H''-E_{R_1}\sim F'\sim D'$ 
and $3\pi''^\ast H'' -\sum_{i=1}^{7}E_{R_i}\sim F''\sim D''$. Similarly $D'+D''$ also gives an closed immersion.

\subparagraph*{} 
Case IV (35 possibilities)

$D'=4\pi^\ast H-\sum_{j=1}^{4}E_{P_{k_j}}-\sum_{i=5}^{7}2E_{P_{l_i}}$ and $D''=3 \pi^\ast H -\sum_{i=1}^{7}E_{P_i}$. In particular 
$D'=4\pi^\ast H-\sum_{j=1}^{4}E_{P_j}-\sum_{i=5}^{7}2E_{P_i}$.
Now, take the quadratic transformation centered at $P_5,P_6,P_7$ and we will get $\widetilde{\mathbb{P}^2}_{P_1P_2...P_7}=\widetilde{\mathbb{P}^2}_{P_1'P_2'P_3'P_4'Q_5Q_6Q_7}$
and $D'\sim F'$ and $D''\sim F''$ where $F'=2\pi'^\ast H'-\sum_{i=1}^4E'_{P_i'}$ and $F''=3 \pi'^\ast H' -\sum_{i=1}^{4}E'_{P_i'}-\sum_{j=5}^{7}E'_{Q_j}$.
Now, we are in the same position as Case II and repeating the argument of Case II, there exist $ R_1,R_2,...,R_7 \in \mathbb{P}^2$ such that $\pi''^\ast H''-E_{R_1}\sim F'\sim D'$ and 
$3\pi''^\ast H'' -\sum_{i=1}^{7}E_{R_i}\sim F''\sim D''$. Then $D'+D''$ also gives an closed immersion.

\subparagraph*{} 
Case V (7 possibilities)

$D'=5\pi^\ast H- E_{P_i}-\sum_{j=1,k_j\neq i}^{6}2E_{P_{k_j}}$ and  $D''=3 \pi^\ast H -\sum_{i=1}^{7}E_{P_i}$. In particular, D=5$\pi^\ast H- E_{P_1}-\sum_{i=2}^{7}2E_{P_i}$.
Now take the quadratic transformation centered at $P_2,P_3,P_4$ and we will get $\widetilde{\mathbb{P}^2}_{P_1P_2...P_7}=\widetilde{\mathbb{P}^2}_{P_1'Q_2Q_3Q_4P_5'P_6'P_7'}$ and 
$D'\sim F'$ and $D''\sim F''$ 
where $F'=4\pi'^\ast H'- E_{P_1'}-\sum_{j=2}^{4}E_{Q_j}-\sum_{i=5}^{7}2E_{P_i'}$ and $F''=3 \pi'^\ast H' -E_{P_1'}-\sum_{i=2}^{4}E'_{Q_i}-\sum_{j=5}^{7}E'_{P_j'}$.
Now we are in the same position as Case IV and repeating the argument of Case IV, there exist $ R_1,R_2,...,R_7 \in \mathbb{P}^2$ such that $\pi''^\ast H''-E_{R_1}\sim F'\sim D'$ and 
$3 \pi''^\ast H'' -\sum_{i=1}^{7}E_{R_i}\sim F''\sim D''$, where
$\pi'':\widetilde{\mathbb{P}^2}_{R_1R_2...R_7}\longrightarrow \mathbb{P}^2$ is the blowing up map at the points  $ R_1,R_2,...,R_7$. Hence $D'+D''$ also gives an closed immersion.

This are all possible very ample divisors of $\widetilde{\mathbb{P}^2}_{P_1P_2...P_7}$ such that this surface can be embedded in $\mathbb{P}^5$ inside the image of $ \mathbb{P}^{1}\times\mathbb{P}^{2}$ given by the Segre embedding and the image of   $p_1\times p_2: \widetilde{\mathbb{P}^2}_{P_1P_2...P_7}\hookrightarrow\mathbb{P}^{1}\times\mathbb{P}^{2}$  is linearly equivalent to (2,2) divisor in $\mathbb{P}^{1}\times\mathbb{P}^{2}$ .
\end{xrem}

\section{Lines of $\mathbb{P}^2$ blown up at seven general points}
We know that $\mathbb{P}^2$ blown up at six general points has 27 lines when we see it as a cubic in $\mathbb{P}^3$ embedded by the anti-canonical divisor.
Also, we know that this lines are all $(-1)$ curves.

Here, we have seen $\mathbb{P}^2$ blown up at seven general points can be embedded in $\mathbb{P}^5$ as well as in $\mathbb{P}^{1}\times\mathbb{P}^{2}$ using the very ample
divisor $D= 4 \pi^\ast H - 2 E_1-\sum_{i=2}^{7}E_i$ or the pair of divisors $(D_1, D_2)$, where $D_1+D_2=D$, and $D_1=\pi^* H-E_1$, $D_2=3\pi^* H-\sum E_i$.  
Note that there are 56 $(-1)$ curves in  $\widetilde{\mathbb{P}^2_7} $; which are

$\bullet$ $L_{ij}=\pi^\ast H-E_i-E_j$, the strict transformation of the line in $\mathbb{P}^2$ containing $P_i$ and $P_j$, $1\leq i, j\leq 7$, (21 possibilities)

$\bullet$ $G_{ij}=\pi^\ast (2H)-\sum_{k\neq i,j} E_k$, the strict transformation of the conic not passing through $P_i$ and $P_j$ 
and passing through the rest of five $P_k$'s of $P_1, P_2,...,P_7$, $1\leq i, j\leq 7$, (21 possibilities)

$\bullet$ $F_i=\pi^\ast (3H)-2 E_i-\sum_{j\neq i} E_j$, the strict transformation of the cubic passing through all the seven points and with a double point at $P_i$
where $i=1,2,...,7$, (7 possibilities).

$\bullet$ Exceptional curves $E_i$, the total transformation of the points $P_i$, $i=1,2,...,7$, (7 possibilities)

\begin{lem}\label{lem6.1}
Any line of $\widetilde{\mathbb{P}^2_7}$ in the embedding of $\mathbb{P}^5$ is a $(-1)$ curve.
 \begin{proof}
Let $L=a \pi^\ast H-\sum_{i=1}^{7}b_iE_i$ be a line in $\widetilde{\mathbb{P}^2_7}$ in the embedding of $\mathbb{P}^5$ where, $a\geq1$, $b_i\geq0$ and $a\geq b_i$ $\forall i$. So $L. (4 \pi^\ast H - 2 E_1-
\sum_{i=2}^{7}E_i)=1 \Rightarrow4a-2b_1-b_2-b_3-b_4-b_5-b_6-b_7=1$ and the genus of $L$ is zero i.e $\frac{(a-1)(a-2)}{2}-\sum_{i=1}^{7}\frac{b_i(b_i-1)}{2}=0$. Solving this two
equations we have, $a^2-\sum_{i=1}^{7} b_i^2=-a+b_1-1$. As $b_1\leq a$, which implies $ -a+b_1-1\leq-1\Rightarrow a^2-\sum_{i=1}^{7} b_i^2\leq-1\Rightarrow L.L\leq-1$.

Let $D$ be any irreducible curve of $\widetilde{\mathbb{P}^2_7}$. The genus of $D$, $g(D)=\frac{1}{2}(D.D-(-K_{S}).D)+1\geq 0$. As the anti-canonical of $\widetilde{\mathbb{P}^2_7}$ 
is an irreducible effective divisor which gives a finite map from $\widetilde{\mathbb{P}^2_7}$ to $\mathbb{P}^2$, $(-K_{S}).D> 0$. Then $D.D\geq -1$. Hence, we have $L.L=-1$.
 \end{proof}

\end{lem}

\begin{xrem}\label{xrem6.2}
In the above, we have listed all 56 $(-1)$ curves of $\widetilde{\mathbb{P}^2_7}$. The divisor gives $L$ is a line in $\widetilde{\mathbb{P}^2_7}$ in the embedding of $\mathbb{P}^5$ if and only if 
$L.(4 \pi^\ast H - 2 E_1-\sum_{i=2}^{7}E_i)=1$. Hence there are only 12 divisor classes in  $\widetilde{\mathbb{P}^2_7}$ which are $E_2, E_3, E_4, E_5, E_6, E_7, L_{12}, L_{13},$ $ L_{14},
L_{15}, L_{16}, L_{17}$ corresponding to lines in the given embedding in $\mathbb{P}^5$. So those are also lines in $\mathbb{P}^{1}\times\mathbb{P}^{2}$. Let $A_1(\mathbb{P}^{1}\times\mathbb{P}
^{2})=\mathbb{Z}\oplus\mathbb{ Z}$ be the group of 1-cycle modulo rational equivalence, which is generated by $pt\times H$ and $\mathbb{P}^1\times pt$. 
Note that any curve rationally equivalent 
to $pt\times H$ or $\mathbb{P}^1\times pt$ is a line in $\mathbb{P}^{1}\times\mathbb{P}^{2}$. This six pairs of lines in $\widetilde{\mathbb{P}^2_7}$ have the property that 
$E_i.L_{1j}=\delta_{ij}$,  $ \forall i, j$, $E_i.E_j=0$ for $i\neq j$ and $L_{1i}.L_{1j}=0$ for $i\neq j$
. We call $(E_i,L_{1i})$ as a pair. So this pair will be either of the form $(\mathbb{P}^1\times pt_1, pt_2\times L_1)$ where $L_1$ is the line passing through $pt_1$ or of the 
form $(pt_3\times L_2, pt_3\times L_3)$ in $\mathbb{P}^{1}\times\mathbb{P}^{2}$. Our claim is that the first situation will never occur. 

Let $pt_2\times L_1 \subset\widetilde{\mathbb{P}^2_7}$ and $pt_2\times\mathbb{P}^2=(1,0)$ be a surface in $\mathbb{P}^{1}\times\mathbb{P}^{2}$.
$(2,2).(1,0).(1,1)=2$ i.e intersection of $\widetilde{\mathbb{P}^2_7}$ with $pt_2\times\mathbb{P}^2$ gives a degree two curve in $\mathbb{P}^5$. Here $pt_2\times L_1\subset
(pt_2\times\mathbb{P}^2)\cap \widetilde{\mathbb{P}^2_7}$ and we already know $pt_2\times L_1$ is a degree one curve. So there is another degree one curve say M which is also inside
$(pt_2\times\mathbb{P}^2)\cap \widetilde{\mathbb{P}^2_7}$.
$M$ is a line in $\mathbb{P}^{1}\times\mathbb{P}^{2}$ and $M\subset pt_2\times \mathbb{P}^2$. So possibilities of $M$ is $pt_2\times L$ where $L$ is a line in $\mathbb{P}^2$.
Clearly $M.(pt_2\times L_1)=1$ in $\widetilde{\mathbb{P}^2_7}$. So $M$ should be another line of the pairing. Hence the claim is proved.
\end{xrem}
\begin{thm}\label{thm6.3}
 As defined in the Remark (\ref{xrem5.3}) , $(\widetilde{\mathbb{ P}^2}_{P_1P_2...P_7}, \mathbb{P}^1, f)$  is a honest conic bundle over $\mathbb{P}^1$.
 \begin{proof}
 To show  $(\widetilde{\mathbb{ P}^2}_{P_1P_2...P_7}, \mathbb{P}^1, f)$  is a honest conic bundle over $\mathbb{P}^1$, we only have to show each fiber is reduced. Assume that, $f$ has
 some non-reduced fiber. Let $F$ be such non reduced fiber. Then $F=2L$ where $L$ is a line in $\mathbb{P}^5$, as degree of $F$ is two in $\mathbb{P}^5$. But $F^2=0$ implies $L^2=0$.
 But this will contradict the Lemma (\ref{lem6.1}). Hence our assumption is not true. So, each fiber of $f$ is either irreducible conic or union of two lines. 
 Hence $(\widetilde{\mathbb{ P}^2}_{P_1P_2...P_7}, \mathbb{P}^1, f)$  is a honest conic bundle over $\mathbb{P}^1$.
 
 \end{proof}

\end{thm}

\begin{corl}\label{corl6.4}
 Lines of $\widetilde{\mathbb{P}^2_7}$ are components of non-irreducible fibers of the conic bundle $(\widetilde{\mathbb{ P}^2}_{P_1P_2...P_7}, \mathbb{P}^1, f)$ over $\mathbb{P}^1$.
 \begin{proof}
  This follows easily from Remark (\ref{xrem6.2}) and Theorem (\ref{thm6.3}).
 \end{proof}

\end{corl}

\section{Smooth surfaces of $\mid(2,2)\mid$ in $\mathbb{P}^{1}\times\mathbb{P}^{2}$}
As $\mid(2,2)\mid$ is a base point free linear system, by Bertini's Theorem, the generic element of the linear system is smooth. The following may be well known,

\begin{thm}\label{thm7.1}
 Any smooth surface of $\mid(2,2)\mid$ of $\mathbb{P}^{1}\times\mathbb{P}^{2}$ can be viewed as the embedding of $\mathbb{P}^2$ blown up at seven  points.
 \begin{proof}
  Let us consider $\mid D\mid=\mid(2,2)\mid$. 

\[\begin{tikzcd}
\mathbb{P}^{1}\times\mathbb{P}^{2} 
\arrow[r, "p_1"] 
\arrow[rd, "p_2"]
& \mathbb{P}^1 \\
& \mathbb{P}^2
\end{tikzcd}\]

where $p_1$, $p_2$ are  two projection maps. $ D=p_1^\ast \mathcal{O}_{\mathbb{P}^1}(2)\otimes p_2^\ast \mathcal{O}_{\mathbb{P}^2}(2)$, 
$-D=p_1^\ast \mathcal{O}_{\mathbb{P}^1}(-2)\otimes p_2^\ast \mathcal{O}_{\mathbb{P}^2}(-2)$. Then\\
$R^i {p_1}_*(p_1^\ast \mathcal{O}_{\mathbb{P}^1}(-2)\otimes p_2^\ast \mathcal{O}_{\mathbb{P}^2}(-2))$\\
= $\mathcal{O}_{\mathbb{P}^1}(-2)\otimes R^i {p_1}_\ast p_2^\ast \mathcal{O}_{\mathbb{P}^2}(-2)$ ( Using projection formula (\cite{H}.III.8))\\
= 0 (as $R^i {p_1}_\ast p_2^\ast \mathcal{O}_{\mathbb{P}^2}(-2)=0$).

Now by the Leray spectral sequence (\cite{H}.III.8)\\
$H^i(\mathbb{P}^{1}\times\mathbb{P}^{2},p_1^\ast \mathcal{O}_{\mathbb{P}^1}(-2)\otimes p_2^\ast \mathcal{O}_{\mathbb{P}^2}(-2))$ \\
$\simeq  H^i(\mathbb {P}^1,{p_1}_\ast(p_1^\ast \mathcal{O}_{\mathbb{P}^1}(-2)\otimes p_2^\ast \mathcal{O}_{\mathbb{P}^2}(-2))$\\
$=H^i(\mathbb{P}^1, \mathcal{O}_{\mathbb{P}^1}(-2)\otimes {p_1}_\ast p_2^\ast \mathscr{O}_{\mathbb{P}^2}(-2))$\\
= 0 for all  $i$, as  ${p_1}_\ast p_2^\ast \mathcal{O}_{\mathbb{P}^2}(-2)=0$

So $h^1(\mathbb{P}^{1}\times\mathbb{P}^{2},\mathcal{O}(-2,-2))=h^2(\mathbb{P}^{1}\times\mathbb{P}^{2},\mathcal{O}(-2,-2))=0$, also we have $h^i(\mathbb{P}^1\times \mathbb{P}^2,\mathcal{O}_{\mathbb{P}^{1}\times\mathbb{P}^{2}})=0$ for $i>0$.

As $D$ is an effective divisor, there is the short exact sequence,
\begin{align}\label{seq4.1}
 0\longrightarrow \mathcal{L}(-D)\longrightarrow \mathcal{O}_{\mathbb{P}^{1}\times\mathbb{P}^{2}}\longrightarrow \mathcal{O}_D\longrightarrow0
 \end{align}
 
Then we get $h^1(\mathcal{O}_D)=h^2(\mathcal{O}_D)=0$
from the induced long exact sequence of cohomologies. Hence the arithmetic genus of the surface $D$, $\rho_a(D)=0$. 
This implies $\chi(\mathcal{O}_D)=1$. We know the canonical divisor of $\mathbb{P}^1\times\mathbb{P}^2$, $K_{\mathbb{P}^1\times\mathbb{P}^2}=  
p_1^\ast K_{\mathbb{P}^1} + p_2^\ast K_{\mathbb{P}^2}= p_1^\ast \mathcal{O}_{\mathbb{P}^1}(-2)+ p_2^\ast \mathcal{O}_{\mathbb{P}^2}(-3)$, where $K_{\mathbb{P}^1}$ and $K_{\mathbb{P}^2}$ are canonical divisor of 
$\mathbb{P}^1$ and $\mathbb{P}^2$ respectively. Note that,\\
$K_D=K_{\mathbb{P}^1\times\mathbb{P}^2} + \mathscr{L}(D)\mid_D$ \textrm{( using the \it adjunction formula\rm)}.\\
$K_D=p_2^\ast \mathcal{O}_{\mathbb{P}^2}(-1)\mid_D$ (as $D\sim (2,2)$).\\
$K_D^2=p_2^\ast \mathcal{O}_{\mathbb{P}^2}(-1)\cdot p_2^\ast \mathcal{O}_{\mathbb{P}^2}(-1)\cdot D= (\mathbb{P}^1 \times pt)\cdot (2 pt\times \mathbb{P}^2+\mathbb{P}^1\times 2H)=2$

By Noether's formula:\\
$12.\chi(\mathcal {O}_D)=\chi_{top}(D)+K_D^2$\\  $\Rightarrow \chi_{top}(D)=10$\\
$\Rightarrow \sum_{i=0}^{4} (-1)^i b_i(D)=10$\\
where $b_i(D)=dim_{\mathbb{ R}} H^i(D,\mathbb{R})$. As $D$ is a surface, $b_0=b_4=1$, and $b_3=b_1$ which deduce 
$$b_2(D)-2b_1(D)=8$$.

We know the irregularity of surface, $ q(D) =h^0(D,\Omega_D)=\frac{1}{2}b_1(D)$, where $\Omega_D$ is the sheaf of differentials on the surface $D$. Also we have the following short exact sequence,
\begin{align}\label{seq2}
0\longrightarrow \mathcal{N}^\ast_{D/{\mathbb{P}^1\times \mathbb{P}^2}}\longrightarrow \Omega_{\mathbb{P}^{1}\times\mathbb{P}^{2}}\otimes\mathcal{O}_D\longrightarrow\Omega_D
\longrightarrow0
\end{align}

where $\mathcal{N}^\ast_{D/{\mathbb{P}^1\times \mathbb{P}^2}}=\mathcal{L}(-D)\mid_D=\mathcal{O}(-2,-2)\mid_D$ is the conormal sheaf of $D$ in $\mathbb{P}^1\times\mathbb{P}^2$.

Note that, $\Omega_{\mathbb{P}^{1}\times\mathbb{P}^{2}}\cong (p_1^\ast \Omega_{\mathbb{P}^1}\oplus p_2^\ast \Omega_{\mathbb{P}^2})$, hence
$h^0(\Omega_{\mathbb{P}^{1}\times\mathbb{P}^{2}})=h^0(p_1^\ast \Omega_{\mathbb{P}^1})+h^0
(p_2^\ast \Omega_{\mathbb{P}^2})=0$.
It can be calculated easily from the long exact sequence of the short exact sequence \ref{seq4.1}, that $h^1(\mathcal{O}(-2,-2)\vert_D)=0$.  
Finally, we get $h^0(D,\Omega_D)=0$ form the long exact sequence of cohomologies induced from (\ref{seq2}). Hence the irregularity of the surface $D$, $q(D)=b_1(D)=0$. Therefore, $b_2(D)=8$.

The second plurigenera of $D$,
$P_2=h^0(D,\mathcal{O}_D(2K_D))=h^0(p_2^\ast \mathcal{O}_{\mathbb{P}^2}(-2)\vert_D)=0$.

So by the Castelnuovo's Rationality Criterion, \cite{BE}, $D$ is a rational surface. We know that every rational surface is either a blow-up of 
$\mathbb{P}^2$ or ruled surface over 
$\mathbb{P}^1$ up to isomorphism (\cite{GH}, page no. 520). In Theorem \ref{thm4.5.4} we observe that any smooth curve of $D$ has self-intersection at least
$-2$. So possibilities of $D$ are blown-up of $\mathbb{P}^2$, $\mathbb{F}_0$, $\mathbb{F}_1$, and $\mathbb{F}_2$, where 
$\mathbb{F}_n=\mathbb{P}(\mathcal{O}_{\mathbb{P}^1}\oplus\mathcal{O}_{\mathbb{P}^1}(-n))$. We know that $\mathbb{F}_1$ is isomorphic to $\mathbb{P}^2$ blow up at one point. Blow up of 
a particular point of $\mathbb{F}_1$ is isomorphic to blow up of a point of $\mathbb{F}_2$ (\cite{GH}, page no. 520). So $\mathbb{F}_2$ blown up at one point is isomorphic to 
$\mathbb{P}^2$ blown up at two points.

Now consider the short exact sequence, 
\begin{align}
 0\longrightarrow \mathbb{Z} \longrightarrow \mathcal{O}_D\longrightarrow \mathcal{O}_D^\ast \longrightarrow0
\end{align}
$Pic(D)\cong H^2(D,\mathbb{Z})$ as we have seen $h^2(\mathcal{O}_D)=0$.
So rank$(Pic(D))=$ rank$H^2(D,\mathbb{Z})=$ dim$_{\mathbb{R}}(H^2(D,\mathbb{R}))= b_2(D)=8$.  We know Picard group of a ruled surface over $\mathbb{P}^1$ is always isomorphic to
$\mathbb{Z}\oplus\mathbb{Z}$. Hence our surface is a blow up 
of $\mathbb{P}^2$ when $Pic(D)=\mathbb{Z}^8$. So $D$ is isomorphic to $\mathbb{P}^2$ blown up at seven points and $$Pic(D)\cong Pic(\mathbb{P}^2)\oplus_{i=1}^{7}\mathbb{Z}.\mathcal{O}( E_i),$$ where 
$E_i$'s are exceptional curves.
 \end{proof}

\end{thm}

\begin{xrem}
 Note that, we have proved that any smooth surface of $\mid(2,2)\mid$ is $\mathbb{P}^2$ blown up at seven points, but they may not be in  general position.
\end{xrem}
\begin{exm}\label{exm7.3}
Let $S$ be a smooth surface in $\mathbb{P}^{1}\times\mathbb{P}^{2}$ from the linear system $\mid(2,2)\mid$, having the equation
$y_0^2(x_0x_1-x_2^2)+y_0y_1(x_1x_2-x_0 ^2)+y_1^2(x_0x_1-x_2^2)$ where $y_0,y_1$ are homogeneous coordinate of $\mathbb{P}^1$, $x_0,x_1,x_2$ are homogeneous coordinate of $\mathbb{P}^2$.
Clearly, $\mathbb{P}^1\times [1,1,1]$ is a line in $\mathbb{P}^{1}\times\mathbb{P}^{2}$ also inside $S$. But in the Corollary
(\ref{corl6.4}), we have seen lines of $\mathbb{P}^2$ blown up at seven general points  are of the form $pt\times L$, where $L$ is a line in $\mathbb{P}^2$. Therefore $S$ is 
$\mathbb{P}^2$ blown up at seven points, which are not in general position.
\end{exm}

\subparagraph*{} Let $D$ be any irreducible divisor of a smooth surface $S$ and $S\sim (2,2)$ in $\mathbb{P}^{1}\times\mathbb{P}^{2}$. The genus of $D$, 
$g(D)=\frac{1}{2}(D.D-(-K_{S}).D)+1\geq 0$. As the anti-canonical is an irreducible effective divisor of $S$, $(-K_{S}).D\geq 0$. Then $D.D\geq -2$. 
If $D.D=-2$, then $D.K_{S}=0$ and $g(D)=0$ i.e $D$ is isomorphic to $\mathbb{P}^1$. We have $S\simeq\widetilde{\mathbb{P}^2_7}$.
If those seven points are not in general position, then $S$ may have some $-2$ lines.

\begin{thm}
 Any smooth surface $S$ of $\mid(2,2)\mid$ of $\mathbb{P}^{1}\times\mathbb{P}^{2}$ has at most four $(-2)$ curves.
 \begin{proof}
  The second projection $p_2$ restricted to $S$ is a generically finite degree two map
   which is defined by the anti-canonical. 
Let $L$ be a $(-2)$ curve in $S$, then $L.(-K_{S})=0$. So the line $L$ is contracted to a point by the morphism $p_2$. Then $L$ will be of the form 
$\mathbb{P}^1\times P_1$ inside 
$\mathbb{P}^{1}\times\mathbb{P}^{2}$. As $S$ is a (2,2) 
surface of $\mathbb{P}^{1}\times\mathbb{P}^{2}$, the defining equation of $S$ is $y_0^2 F_0(x_0,x_1,x_2)+y_1^2 F_1(x_0,x_1,x_2)+ y_0y_1 F_2(x_0,x_1,x_2)$,
where $y_0,y_1$
are homogeneous co-ordinate of $\mathbb{P}^1$, $x_0,x_1,x_2$ are homogeneous co-ordinate of $\mathbb{P}^2$ and deg$(F_i(x_0,x_1,x_2))=2$, for $i=0,1,2$.
If $\mathbb{P}^1\times P_1\subset
Z(y_0^2 F_0(x_0,x_1,x_2)+y_1^2 F_1(x_0,x_1,x_2)+ y_0y_1 F_2(x_0,x_1,x_2))$, then $P_1\in Z(F_0,F_1,F_2)$.

Conversely, if $P_1\in Z(F_0,F_1,F_2)$, then $L=\mathbb{P}^1\times P_1 \subset S$ and $L$ is contracted to $P_1$ by $p_2$ i.e $L.(-K_{S})=0$, then $L.L=-2$. Any $-2$ curve  of $S$, 
which is always a line, will be of the form $\mathbb{P}^1\times P$ where $P\in Z(F_0,F_1,F_2)$. But we know $\#Z(F_0,F_1,F_2)\leq4$, as $F_i$ are degree two curves in $\mathbb{P}^2$.
Hence the result is proved.
 \end{proof}

\end{thm}

\begin{xrem}
  $(-2)$ lines of $S$ are reduced and in the form of $\mathbb{P}^1\times P_i$. Hence they are disjoint to each other.
\end{xrem}
\begin{exm}
  $Z(y_0^2(x_0x_2-x_1^2)+y_1^2(x_0x_1-x_2^2))$ is a smooth surface which is linearly equivalent to (2,2) of $\mathbb{P}^{1}\times\mathbb{P}^{2}$. It has four $(-2)$ lines, which are 
$\mathbb{P}^1\times [1,1,1]$, $\mathbb{P}^1\times [1,0,0]$, $\mathbb{P}^1\times [1,\omega,\omega^2]$, $\mathbb{P}^1\times [1,\omega^2,\omega]$, where $\omega$ is a cubic root of unity.
\end{exm}
\begin{prop}
 Suppose that a smooth surface $S$ of $\mid(2,2)\mid$ of $\mathbb{P}^{1}\times\mathbb{P}^{2}$ has a $(-2)$ curve $L$. Then either $L\sim E_i-E_j$ or 
$L\sim \pi^\ast H-E_i-E_j-E_k $ or $L\sim 2\pi^\ast H-\sum_{i=1}^{6} E_{k_i}$ 
\begin{proof}
Let $L$ be a $(-2)$ curve, then either $L$ is component of some exceptional curve or $L=a\pi^\ast H-\sum_{i=1}^7 b_iE_i$, where $a\geq 1$ or $b_i\geq0$.
If $L$ is a component of an exceptional curve, then the only possibility for $L$ is $E_i-E_j$. Otherwise $L^2=a^2-\sum_{i=1}^7 b_i^2=-2$ and $L.K_S=0$ implies $3a=\sum_{i=1}^7b_i$.

By Schwarz's inequality $\mid\sum_{i}x_iy_i\mid \leq\mid\sum_{i}x_i^2\mid .\mid\sum_{i}y_i^2\mid$

Here $x_i=1, y_i=b_i$ for i=1,2,...7 and $x_i=0, y_i=0$ for i$>$7. So we have

$\Rightarrow(\sum b_i)^2 \leq 7(\sum b_i^2)$

$\Rightarrow (3a)^2\leq 7(a^2+2)$

$\Rightarrow a^2\leq 7$

$\Rightarrow a\leq 2$

Then by the above equations, either $L\sim \pi^\ast H-E_i-E_j-E_k $ or $L\sim 2\pi^\ast H-\sum_{i=1}^{6} E_{k_i}$. 
\end{proof}
\end{prop}

\subparagraph*{}
In the fifth section, we have seen all possible very ample divisors, by which $\mathbb{P}^2$ blown up at seven general points can be embedded as (2,2) divisor in 
$\mathbb{P}^{1}\times\mathbb{P}^{2}$ as well as in $\mathbb{P}^5$. In this section we already have seen that there are other smooth surfaces linearly equivalent to (2,2) divisor,
which again can be viewed as $\mathbb{P}^2$ blown up seven points, where those seven points are not in general position. Now, we will see, what is a very ample divisors in this cases.
(See Remark \ref{xrem7.12}).

\subparagraph*{}

Consider the surface $S$ for which we have a sequence of morphism 

$S=S_n\xrightarrow{\pi_n} S_{n-1}\xrightarrow{\pi_{n-1}} S_{n-2}....\xrightarrow{\pi_2} S_{1}\xrightarrow{\pi_1} S_{0}=\mathbb{P}^2$

where $S_i\xrightarrow{\pi_i} S_{i-1}$ is the blowing-up of $S_{i-1}$ at a point $P_i$. $\mathscr{E}=\{\mathscr{E}_0,\mathscr{E}_1,...,\mathscr{E}_n\}$ can be considered as a free basis
of Pic$S$, where $\mathscr{E}_0$ is the class of a line in $\mathbb{P}^2$ and $\mathscr{E}_i$ is the class of $E_i=\pi_i^{-1}(P_i)$. Such collection $\mathscr{E} $ of divisor classes
is called an \textit{ exceptional configuration}. Thus, there is a bijection between a sequence of morphism given above from $S$ to $\mathbb{P}^2$ and an exceptional configuration of $S$.

\begin{thm}\label{thm7.9}
 Let $\mathscr{L}\in$ Pic$S$. Then $\mathscr{L}$ is very ample iff there is an exceptional configuration $\mathscr{E}=\{\mathscr{E}_0,\mathscr{E}_1,...,
\mathscr{E}_n\}$ of $S$ such that (i) $\mathscr{L}.(\mathscr{E}_0-\mathscr{E}_1)>0$, (ii) $\mathscr{L}.(\mathscr{E}_0-\mathscr{E}_1-\mathscr{E}_2)>0$, (iii)
$\mathscr{L}.(\mathscr{E}_0-\mathscr{E}_1-\mathscr{E}_2-\mathscr{E}_3)\geq0$, (iv) $\mathscr{L}.(\mathscr{E}_i-\mathscr{E}_{i+1})\geq0$ for $i\geq1$, (v) $\mathscr{L}.\mathscr{E}_i>0$, 
(vi) $\mathscr{L}.\mathscr{N}>0$ for any $-2$ curve $\mathscr{N}$, (vii) $\mathscr{L}.K_S\leq-3$ where $\mathscr{E}_0=\pi^\ast H$.
\begin{proof}
 \cite{Ha} Theorem (2.1).
\end{proof}
\end{thm}

\begin{xrem}\label{xrem7.10}
 Let $S$ and $S'$ be two surfaces, where $S$ is $\mathbb{P}^2$ blown up at $n$ general points and $S'$ is $\mathbb{P}^2$ blown up at $n$ points which are not in general position. Let 
$\mathscr{L}'$ be a very ample divisor of $S'$. By the Theorem (\ref{thm7.9}) there exist an exceptional configuration $\mathscr{E}'$ such that $\mathscr{L}'$ satisfies (i)-(vii)
and let $\mathscr{L}'=a\mathscr{E}_0'-\sum_{i=1}^7 b_i\mathscr{E}_i'$. Now, let $\mathscr{L}=a\mathscr{E}_0-\sum_{i=1}^7 b_i\mathscr{E}_i$ be a divisor of $S$ with
respect to the exceptional configuration $\mathscr{E}$, then $\mathscr{L}$ satisfies all the properties of Theorem (\ref{thm7.9}), which implies that $\mathscr{L}$ is also a very ample divisor 
of $S$.

The Picard group of $\widetilde{\mathbb{P}^2_r} $ is $\mathbb{Z}l\oplus_{i=1} ^r \mathbb{Z} E_i$, where $l$ is $\pi^\ast (H)$. Let $D= a \pi^\ast H -\sum_{i=1}^{7} b_i E_i$ be a divisor
of $\widetilde{\mathbb{P}^2_r} $ which is also denoted as $(a, b_1, b_2,...,b_r)$ as an element of $\mathbb{Z}^{r+1}$. So by the above argument,
if $(a, b_1, b_2,...,b_r)$ is a very ample divisor of $\widetilde{\mathbb{P}}^2_{P_1P_2..P_r} $ for some given set of $r$ points in $\mathbb{P}^2$, 
then $(a, b_1, b_2,...,b_r)$ is also very ample divisor of the surface $\mathbb{P}^2$ blown-up at $r$ general points. But the converse (if $(a, b_1, b_2,...,b_r)$ is a very ample 
divisor for $\widetilde{\mathbb{P}}^2_{P_1P_2..P_r} $ where $r$ points are in general position in $\mathbb{P}^2$ whether $(a, b_1, b_2,...,b_r)$ is  a very ample divisor
of $\widetilde{\mathbb{P}}^2_{P_1P_2..P_r} $ for some given set of $r$ points in $\mathbb{P}^2$ or not) is not true in general because of the property (vi) of Theorem (\ref{thm7.9}). 
\end{xrem}

\begin{xrem}\label{xrem7.11}
 In the Remark (\ref{xrem5.4}), we listed all possible very ample divisors of $\mathbb{P}^2$ blown up at seven general points which give closed immersion of the surface in 
$\mathbb{P}^{1}\times\mathbb{P}^{2}$ as (2,2) type divisor. Also, in the Remark (\ref{xrem5.4}) we noted that, if $\mathscr{L}$ is such a very ample divisor then there exist a exceptional
configuration $\mathscr{E}=\{\mathscr{E}_0,\mathscr{E}_1,...,\mathscr{E}_n\}$ such that  $\mathscr{L}$ can be written as $4\mathscr{E}_0-2\mathscr{E_1}-\sum_{i=2}^7\mathscr{E}_i$. 

\end{xrem}

\begin{xrem}\label{xrem7.12}
 Let $S$ be a smooth surface $S\sim (2,2)$ of $\mathbb{P}^{1}\times\mathbb{P}^{2}$ and $S\simeq \widetilde{\mathbb{P}}^2_{P_1P_2...P_7} $, where $P_1,P_2,...,P_7$ are not in 
general position. Also, let $\mathscr{L}$ be the very ample divisor of $\widetilde{\mathbb{P}}^2_{P_1P_2...P_7} $ which gives the above closed immersion. Then by Remarks (\ref{xrem7.10}) and
(\ref{xrem7.11}) there is an exceptional configuration $\mathscr{E}'=\{\mathscr{E}_0',\mathscr{E}_1',...,\mathscr{E}_n'\}$ of $\widetilde{\mathbb{P}}^2_{P_1P_2...P_7} $ such that
$\mathscr{L}$ can be written as $4\mathscr{E}_0'-2\mathscr{E_1}'-\sum_{i=2}^7\mathscr{E}_i'$. 
\end{xrem}

\begin{exm}
 Let us consider seven points from $\mathbb{P}^2$ such that $P_2, P_3, P_4$ are collinear, $P_4, P_5, P_6$ are collinear, $P_2, P_7, P_5$ are collinear,
$P_3, P_6, P_7$ are collinear and $P_1$ is not collinear with any two $P_i$'s and no six points lie on a conic. Let $S$ be a surface obtained from blown up of $\mathbb{P}^2$ at $P_1, P_2,..
.,P_7$. Here, $\pi^\ast H-E_2-E_3-E_4$, $\pi^\ast H-E_4-E_5-E_6$, $\pi^\ast H-E_2-E_7-E_5$ and 
$\pi^\ast H-E_3-E_7-E_6$ are $-2$ lines of $S$. Then,  $ 4 \pi^\ast H - 2 E_1-\sum_{i=2}^{7}E_i$ is a very ample divisor of $S$ (Theorem \ref{thm7.9}).
Clearly, $ 3 \pi^\ast H -\sum_{i=1}^{7}E_i$ gives generic degree two map from $S$ to $\mathbb{P}^2$ and $\pi^\ast H-E_1$ gives conic bundle
map from $S$ to $\mathbb{P}^1$. Hence, $S$ is a smooth (2,2) surface of $\mathbb{P}^1\times\mathbb{P}^2$ which has four $-2$ lines.

\end{exm}

\end{document}